\documentclass[a4paper, 12pt]{amsart}
\usepackage{fourier}
\usepackage{amssymb}
\usepackage{amsmath}
\usepackage{amscd}
\usepackage{amsthm}
\usepackage[centertags]{amsmath}
\usepackage{amsfonts}
\usepackage{newlfont}
\usepackage[all]{xy}
\usepackage{graphicx}
\usepackage{amsfonts, amssymb}
\usepackage[usenames]{color}
\usepackage{mathrsfs}
\usepackage{latexsym}

\newtheorem{thm}{Theorem}[section]
\newtheorem{cor}[thm]{Corollary}
\newtheorem{lem}[thm]{Lemma}
\newtheorem{prop}[thm]{Proposition}

\theoremstyle{definition}

\newtheorem{rem}[thm]{Remark}
\newtheorem{exa}[thm]{Example}
\newtheorem{notation}[thm]{Notation}

\newcommand{\zC}{\mathbb C}

\newcommand{\zR}{\mathbb R}
\newcommand{\zN}{\mathbb N}
\newcommand{\zK}{\mathbb K}

\begin{document}

\title[On the rank and the approximation of symmetric tensors]{On the rank and the approximation of symmetric tensors}

\thanks{This work was partially supported by CONICET PIP 11220130100329  and ANPCyT PICT 2018-04250.} 
\thanks{Formal publication  at https://doi.org/10.1016/j.laa.2021.07.002}


\author[J. T. Rodr\'{i}guez]{Jorge Tom\'as Rodr\'{i}guez}
\address{Departamento de Matem\'{a}tica and NUCOMPA, Facultad de Cs. Exactas, Universidad Nacional del Centro de la Provincia de Buenos Aires, (7000) Tandil, Argentina and CONICET}
\email{jtrodrig@dm.uba.ar}

\begin{abstract}
In this work we study different notions of ranks and  approximation  of tensors. We consider the  tensor rank, the nuclear rank and we introduce the notion of symmetric decomposable rank, a notion of rank defined only on symmetric tensors. We show that when approximating symmetric tensors, using the symmetric decomposable rank has some significant advantages over the tensor rank and the nuclear rank. 
\end{abstract}

\maketitle

\section{Introduction}
For $\zK$, the field of complex numbers $\zC$ or the field of real numbers $\zR$, the space of multiway arrays (also known as multidimensional arrays) $\zK^{n_1\times\cdots\times n_d}$ fulfils an important role in several application areas such as Chemometrics, Signal Processing, Data Compression, or Data Analysis.  In several of these applications, it is important to approximate multiway arrays  by  simpler multiway arrays that are sparsely representable. By identifying the space $\zK^{n_1\times\cdots\times n_d}$ with the tensor product $\bigotimes^d \zK^{n_i}$ this problem translates into approximating tensors  with low-rank tensors. Therefore, given a tensor $\mathbf{z} \in \bigotimes^d \zK^{n_i}$,  sometimes is of use to find a rank$-r$ approximation of $\mathbf{z}$. That is, finding a tensor $\mathbf{x}$ of rank at most $r$ such that $$\alpha(\mathbf{z} -\mathbf{x})$$ is small, where $\alpha(\, \cdot \,)$ is a norm on $\bigotimes^d \zK^{n_i}$ and rank is a notion that in some sense measures how complex is a tensor. Two reasons to approximate tensors are the same as when working with matrices. In general, low-rank tensors are more inexpensive to store, and also, making computations with low-rank tensors is faster. For a more detailed explanation on the need for low-rank approximations, we refer the reader to  \cite{cichocki2015tensor, comon2009tensor, friedland2015low, grasedyck2013literature, sidiropoulos2017tensor}, and the references therein.

It is also of interest to know if there is a best  rank$-r$ approximation of $\mathbf{z}$: a tensor $\mathbf{x}$ of rank at most $r$ such that
$$\alpha( \mathbf{z}-\mathbf{x}) = \inf \{\alpha( \mathbf{z}-\mathbf{y}): \mathbf{y}\in \textstyle{\bigotimes^d}\zK^{n_i}, \mathbf{y} \text{ has rank at most } r\}.$$
Naturally, this is a necessary condition for algorithms designed to find a best  rank$-r$ approximation to work. But in some papers, the existence of such a tensor is wrongly assumed.  This lack of rigor has been pointed out in  \cite{comon2008symmetric, de2008tensor}.

The trade-off done when a tensor $\mathbf{z}$ is replaced by an approximation $\mathbf{x}$ is that some information may be lost. Therefore, it can be useful to develop results to recover as much information as possible of $\mathbf{z}$ using only its approximation $\mathbf{x}$.

The present work aims to investigate the problems exposed above, taking into account different norms. Rather than creating algorithms, we focus more on developing the theoretical frame needed by these algorithms. Given a tensor $\mathbf{z} \in \bigotimes^d \zK^{n_i}$, its tensor rank is defined as
$$
\operatorname{rank}(\mathbf{z}):= \min\left\{r\in \zN: \mathbf{z}=\sum_{i=1}^r z_1^i\otimes \cdots \otimes z_d^i\right\}.
$$
For this notion of rank, we study the best rank$-1$ approximation of a symmetric tensor $\mathbf{z}$ considering the Hilbert-Schmidt norm. In Propositions \ref{prop best rank 1} and \ref{prop best rank 1 comp} we show that if \mbox{$x_1\otimes \cdots \otimes x_d$} is a best rank$-1$ approximation of $\mathbf{z}$ then the vectors $x_1,\ldots, x_n$ are either collinear or coplanar in the real case, while they are necessarily collinear in the complex case. In particular, in the complex case, a best rank$-1$ approximation of a symmetric tensor $\mathbf{z}$ has to be symmetric.

In the real case, when the vectors $x_1,\ldots, x_n$ are coplanar, in Theorem \ref{teo recover} we provide an algorithm to obtain a best rank$-1$ approximation of $\mathbf{z}$ of the form
$$v\otimes v\otimes w\otimes \cdots \otimes w,$$
with $v$ and $w$ orthonormal. This allows to partially recover the tensor $\mathbf{z}$ using only its best rank$-1$ approximation. Part of this algorithm is obtained following the ideas in the proof of \cite[Theorem 2.1]{pappas2007norm}. This part can also be used to get a best rank$-1$ approximation of $\mathbf{z}$ of the form
$$x\otimes  \cdots \otimes x.$$
The main importance of such an approximation is that it is much cheaper to store than the original one, the only information needed is the vector $x$.

A downside of the tensor rank is that there is not always a best rank-$r$ approximation for $r>1$. A notion of rank that does not have this problem  is the nuclear rank introduced in \cite{friedland2014nuclear}. This notion is closely related with the projective norm. Given a tensor $\mathbf{z} \in \bigotimes^d \zK^{n_i}$ its nuclear rank is defined as
$$\operatorname{rank}_\pi(\mathbf{z}):= \min\left\{r\in \zN: \mathbf{z}=\sum_{i=1}^r z_1^i\otimes \cdots \otimes z_d^i, \pi(\mathbf{z})=\sum_{i=1}^r \Vert z_1^i \Vert \cdots \Vert z_d^i\Vert\right\},$$
where 
$$\pi(\mathbf{z}):= \inf \left\{\sum_{i=1}^r  \Vert z_1^i\Vert \cdots \Vert z_d^i\Vert: \mathbf{z} = \sum_{i=1}^r   z_1^i\otimes \cdots \otimes z_d^i \right\}$$
is the usual projective norm. We investigate the representations that give the projective norm, and therefore the ones needed to compute the nuclear rank. Namely, we show that if $\mathbf{z}$ is symmetric, and
$$\mathbf{z}=\sum_{i=1}^r z_1^i\otimes \cdots \otimes z_d^i\,\,\,\text{ with }\,\,\, \pi(\mathbf{z})=\sum_{i=1}^r \Vert z_1^i \Vert \cdots \Vert z_d^i\Vert,$$
then for each $i=1,\ldots, n$, the vectors $z_1^i,\ldots,z_d^i$ are either collinear or coplanar in the real case, while they are necessarily collinear in the complex case. This significantly reduces the possible representations of $\mathbf{z}$ that need to be considered to compute its projective norm.

In the final part of the present work, we introduce a new notion of rank, that we call decomposable symmetric rank. For  a symmetric tensor $\mathbf{z}~\in~\bigotimes^{d,s}~\zK^n$, we define its decomposable symmetric rank by
$$\operatorname{rank}_\sigma(\mathbf{z}):= \min\left\{r\in \zN: \mathbf{z}=\sum_{i=1}^r z_1^i\vee \cdots \vee z_d^i\right\},$$
where $z_1^i\vee \cdots \vee z_d^i$ is the  decomposable symmetric tensor obtained by summing over the group of permutations of $\{1,\ldots,d\}$
$$z_1^i\vee \cdots \vee z_d^i:= \frac{1}{d!}\sum_{\eta \in S_d} z^i_{\eta(1)}\otimes \cdots \otimes z^i_{\eta(d)}.$$
This notion is similar to the tensor rank, the difference is that instead of considering sums of elementary tensors, we consider sums of decomposable symmetric tensors. When we want to approximate a symmetric tensor, this notion of rank has some advantages over the tensor rank and the nuclear rank. In Theorem \ref{teo aprox dec} we show that if $\mathbf{z}$ is a symmetric tensor, and $\mathbf{y}$ is a tensor rank$-r$ approximation of $\mathbf{z}$, there is an easy way to construct a decomposable symmetric rank$-l$ tensor $\mathbf{x}$ such that $l\leq r$ and that $\mathbf{x}$ is an approximation at least as good as $\mathbf{y}$. This holds true for either the Hilbert-Schmidt, the injective, and the projective norms. The fact that $l\leq r$ implies that storing $\mathbf{x}$ is no more expensive than storing $\mathbf{y}$. In Example \ref{exam rank 2} we show that, as happens for the tensor rank, there is not always a best rank$-r$ approximation for the decomposable symmetric rank. But this is not the case for $r=1$. As a consequence of Lemma \ref{best rank 1 equivalence dec}, there is always a best decomposable symmetric rank$-1$ approximation.

The article is organized as follows. In Section \ref{sec prel} we fix some notation and do a quick overview of the theory of tensor products needed to develop our results. In Section \ref{sec tensor rank} we  state and prove our results regarding the tensor rank. The results on the nuclear rank can be found in Section \ref{sec nuclear rank}. Finally, in Section \ref{sect dec rank} we introduce  the decomposable symmetric rank and give the results concerning this notion.

\section{Preliminaries}\label{sec prel}

Before venturing into the main topic of study, let us fix some notation and definitions related to the tensor product and the symmetric tensor product. In this section, we also review some of the theory regarding these topics that will be needed  in the development of this work.

\subsection*{Tensor product}

Let $\zK$ be either the field of complex numbers $\zC$ or the field of real numbers $\zR$ and $n_1,\ldots, n_d$ positive integers. We study the \emph{tensor product}
$$\textstyle{\bigotimes^d} \zK^{n_i}.$$
When there is $n$ such that $n_i=n$ for every $i$, we simply  write 
$$\textstyle{\bigotimes^{d}} \zK^n.$$

For elements of $\zK^{n_i}$ we use plain letters $x, y, z$, while for tensors we use bold letters $\mathbf{x}, \mathbf{y}, \mathbf{z}$.  A tensor $\mathbf{x}$ is called \emph{elementary tensor} if there are vectors $x_1\in \zR^{n_1},\ldots, x_d \in \zR^{n_d}$ such that 
$$\mathbf{x}=x_1\otimes \cdots \otimes x_d.$$
In the literature, these tensors are often called decomposable. We will not use this terminology since we are going to use a  similar one for symmetric tensors, and we want to avoid confusion.

For a modern exposition on the tensor product of normed vector spaces, we refer the reader to \cite{defant1992tensor, diestel2008metric, ryan2002introduction}.

\subsection*{Symmetric tensor product}

We write $$\textstyle{\bigotimes^{d,s} \zK^n}$$ for the \emph{symmetric tensor product}, the subspace of $\bigotimes^d \zK^{n}$ consisting on all the symmetric tensors. 

The \emph{symmetrization operator} 
$$\textstyle{\sigma:\bigotimes^d \zK^n \rightarrow \bigotimes^{d,s} \zK^n}$$
is useful to relate the tensor product with the subspace of symmetric tensors. This operator is defined on elementary tensors by
$$\sigma(x_1\otimes \cdots \otimes x_d) =x_1\vee \cdots \vee x_d := \frac{1}{d!}\sum_{\eta \in S_d} x_{\eta(1)}\otimes \cdots \otimes x_{\eta(d)},$$
where $S_d$ is the group of permutations of $\{1,\ldots,d\}$. A symmetric tensor $\mathbf{x}$ is called \emph{decomposable symmetric tensor} if there are vectors $x_1,\ldots, x_d\in \zK^n$ such that
$$\mathbf{x}=x_1\vee \cdots \vee x_d.$$
That is, the decomposable symmetric tensors are the image of the elementary tensors via the symmetrization operator $\sigma$.

For a deeper introduction to symmetric tensor products, we refer the reader to Floret’s survey \cite{floret1997natural}.

\subsection*{Relation with multilinear forms and polynomials}

There is a natural way of identifying $\zK^n$ with its dual $$(\zK^n)^*=\{f:\zK^n\rightarrow \zK: f \text{ is linear}\}.$$ This is obtained by identifying an element $x\in \zK^n$ with the linear function $f_x$ given by 
$$f_x(y):=\langle y, x \rangle,$$
where $\langle  x  , y  \rangle=y^*x$ is the usual inner product on $\zK^n$. 

Similarly,  we can find a bijection between the tensor product $\bigotimes^d \zK^{n_i}$ and the space of multilinear forms $$\mathcal{L}(\zK^{n_1}\times\cdots\times \zK^{n_d})=\{\operatorname L:\zK^{n_1}\times\cdots\times \zK^{n_d}\rightarrow \zK:\operatorname L \text{ is multilinear}\}.$$ In order to do this, a tensor $\mathbf{z} = \sum_{i=1}^r   z_1^i\otimes \cdots \otimes z_d^i$ is associated with the multilinear form $\operatorname{L}_\mathbf{z}$ defined as
\begin{equation}\label{ec ident mult}
\operatorname{L}_\mathbf{z}(y_1,\ldots,y_d):=\sum_{i=1}^r   \langle y_1, z_1^i\rangle \cdots \langle y_d, z_d^i\rangle.
\end{equation}

This procedure  can also be used to identify the symmetric tensor product $\bigotimes^{d,s} \zK^{n}$ with the space of $d$-homogeneous polynomials $$\mathcal{P}(^d\zK^{n})=\{\operatorname P:\zK^n\rightarrow \zK:\operatorname P \text{ is a $d$-homogeneous polynomial}\}.$$ A symmetric tensor $\mathbf{u}=\sum_{i=1}^r   u^i\otimes \cdots \otimes u^i$ is associated to the polynomial
\begin{equation}\label{ec ident pol}
\operatorname{P}_\mathbf{u}(y):=\sum_{i=1}^r   \langle y, u^i\rangle^d.
\end{equation}

We remark that in the complex case these bijections are anti linear. It is not hard to change them slightly to obtain linear identifications. For the purposes of this paper either can be used. We chose these since they are lighter on the notation when we write everything in terms of inner products.

\subsection*{Polarization constant} The \emph{$d$-th polarization constant} of a normed space $E$, noted as $c(d,E)$, is defined as the infimum of all the constants $C$ such that 
\begin{equation}\label{def polarization}
\sup \{|\operatorname L(y_1,\ldots,y_d)|:\Vert y_i\Vert = 1 \forall i\}\leq C \sup \{|\operatorname L(y,\ldots,y)|: \Vert y \Vert =1\}
\end{equation}
for any continuous multilinear symmetric form $\operatorname L:\underbrace{E\times \cdots \times E}_{d\,\,\, times}\rightarrow \zK$. 

This constant is closely related to the metric theory of tensor products. We will not work directly with this constant for an arbitrary normed space, and its definition is given with an informative purpose only. However, for some results, we will rely on the fact that in the particular case of Hilbert spaces the polarization constant is one. This was proved by Banach \cite{banach1938homogene} (for further references see  \cite[Section 2.1]{floret1997natural}).

\subsection*{Norms in tensors products} 

Next we recall some basic norms on tensor products. Along this article the norm considered on $\zK^{n_i}$ is always the euclidean norm. On $\bigotimes^d \zK^{n_i}$ define the inner product on elementary tensors by
\begin{equation}\label{inner product}
\langle y_1 \otimes \cdots \otimes y_d, x_1 \otimes \cdots \otimes x_d \rangle = \prod_{i=1}^d \langle y_i, x_i \rangle,
\end{equation}
where $\langle y_i, x_i \rangle=y_i^*x_i$. This inner product gives the \emph{Hilbert-Schmidt norm} for tensors
$$\operatorname{HS}(\mathbf{z}):= \sqrt{ \langle \mathbf{z}, \mathbf{z} \rangle }.$$

Other norm often used in tensor products is the \emph{injective tensor norm}, defined as
$$\varepsilon(\mathbf{z}) := \sup \{|\langle  y_1\otimes \cdots \otimes y_d, \mathbf{z}\rangle | : \Vert y_i \Vert =1, i =1,\ldots, d\}.$$
By \eqref{ec ident mult}, the injective norm coincides with the spectral norm of $\operatorname{L}_\mathbf{z}$
$$\Vert \operatorname{L}_\mathbf{z}\Vert=\sup\left\{|\operatorname{L}_\mathbf{z}(y_1,\ldots, y_d)|:\Vert y_i\Vert = 1, i =1,\ldots, d\right\}.$$ 
For this reason, the injective norm  is sometimes also referred to as the  spectral norm.

The last norm on $\bigotimes^d \zK^{n_i}$ that will appear on this work is the \emph{projective tensor norm}. This norm is given by
$$\pi(\mathbf{z}) := \inf \left\{\sum_{i=1}^r  \Vert z_1^i\Vert \cdots \Vert z_d^i\Vert: \mathbf{z} = \sum_{i=1}^r   z_1^i\otimes \cdots \otimes z_d^i \right\}.$$
In this context, the projective tensor norm of a tensor $\mathbf{z}$ coincides with the nuclear norm of the multilinear form $\operatorname{L}_\mathbf{z}$. Because of this, in the literature $\pi(\, \cdot \,)$ is also often called the nuclear norm.

The injective and projective tensor norms are closely related by the inner product defined in \eqref{inner product}. For two tensors $\mathbf{x},\mathbf{y}\in \bigotimes^d \zK^{n_i}$ we have
\begin{equation}\label{duality}
|\langle \mathbf{y},\mathbf{x}\rangle| \leq \pi(\mathbf{y}) \varepsilon(\mathbf{x}).
\end{equation}
Moreover, the injective and projective tensor norms can be computed as follows
$$\varepsilon(\mathbf{z})=\sup\{ |\langle \mathbf{y},\mathbf{z}\rangle|: \pi(\mathbf{y})=1\}$$
$$\pi(\mathbf{z})=\sup\{ |\langle \mathbf{z},\mathbf{x}\rangle|: \varepsilon(\mathbf{x})=1\}.$$

These norms have their symmetric counterpart. If $\mathbf{z}\in \bigotimes^{d,s} \zK^n$ is a symmetric tensor, its \emph{injective $s$-tensor norm} is defined as
\begin{equation}\label{sim iny}
\varepsilon_s (\mathbf{z}):= \sup \{|\langle  y\otimes \cdots \otimes y, \mathbf{z} \rangle | : \Vert y \Vert =1 \}.
\end{equation}
By \eqref{ec ident pol}, this  norm coincides with he uniform norm of the polynomial $\operatorname{P}_\mathbf{z}$ associated to the tensor
$$\Vert \operatorname{P}_\mathbf{z}\Vert=\sup\{|\operatorname{P}_\mathbf{z}(y)|:\Vert y\Vert = 1\}.$$ 
The \emph{projective $s$-tensor norm} is given by 
\begin{equation}\label{sim proy}
\pi_s(\mathbf{z}) := \inf \left\{\sum_{i=1}^r  \Vert z^i\Vert^d: \mathbf{z} = \sum_{i=1}^r   z^i\otimes \cdots \otimes z^i \right\}.
\end{equation}

A classical result of Banach \cite{banach1938homogene} states that on Hilbert spaces the norm of a symmetric multilinear form coincides with the norm of its associated polynomial (for a modern exposition on this result see \cite{bochnak1971polynomials,pappas2007norm}). In this context, this means that for a symmetric tensor $\mathbf{u},$ the norm of $\operatorname{L}_\mathbf{u}$ coincides with the norm of $\operatorname{P}_\mathbf{u}$. This implies that for Hilbert spaces, the  injective $s$-tensor norm $\varepsilon_s(\, \cdot \,)$ is just the usual injective tensor norm $\varepsilon(\, \cdot \,)$ restricted to symmetric tensors.

As a consequence of the aforementioned result due to Banach, the polarization constant of Hilbert spaces is one. Therefore, the projective $s$-tensor norm $\pi_s(\, \cdot \,)$ also is the restriction of the projective tensor norm $\pi(\, \cdot \,)$ to symmetric tensors (see \cite[Section 2.3]{floret1997natural}). For an alternative proof of this fact, without the direct use of polarization constants, see \cite[Section 5]{friedland2014nuclear}.

Hence, \eqref{sim iny} and \eqref{sim proy} are an alternative way of computing the injective and projective tensor norms for Hilbert spaces.

\section{Tensor rank}\label{sec tensor rank}

For a tensor $\mathbf{z} \in \bigotimes^d \zK^{n_i}$, as mentioned its \emph{tensor rank} is defined as
$$\operatorname{rank}(\mathbf{z}):= \min\left\{r\in \zN: \mathbf{z}=\sum_{i=1}^r z_1^i\otimes \cdots \otimes z_d^i\right\}.$$
This is one of the most common, and probably one of the most important, notions of rank.

A problem with the tensor rank is that there is not always a best rank-$r$ approximation (see \cite{comon2008symmetric, de2008tensor}). However, this is not the case for $r=1$. There is always a best  tensor rank-$1$ approximation. For the Hilbert-Schmidt norm, this problem has been studied by several authors \cite{friedland2013best, friedland2014number, zhang2012best}. 

The following result characterizes the best rank-$1$ approximation. A proof for $\zK=\zR$ can be found in \cite[Lemma 4]{friedland2013best}, but the same technique can be applied to prove the complex case. We use a similar idea later on,  in the proof of Lemma~\ref{best rank 1 equivalence dec}.

\begin{lem}\label{best rank 1 equivalence} Let $\mathbf{z} \in \bigotimes^d \zK^{n}$ and $\mathbf{x}=\lambda x_1\otimes\cdots \otimes x_d,$ with  $x_1, \ldots ,x_d\in \zK^d$ norm one vectors and $\lambda\in \zK$, an elementary tensor. Then $\mathbf{x}$  is a best tensor rank-$1$ approximation of  $\mathbf{z}$ for the Hilbert-Schmidt norm if and only if
$$\lambda=\langle \mathbf{z}, x_1\otimes \cdots \otimes x_d\rangle\,\,\,\text{ and } \,\,\,|\lambda|=\varepsilon(\mathbf{z}).$$
\end{lem}

In terms of multilinear  forms this lemma states that 
$$\overline{\lambda}= \operatorname{L}_\mathbf{z}(x_1,\ldots, x_d)\,\,\,\text{ and } \,\,\,|\lambda|=\Vert \operatorname{L}_\mathbf{z}\Vert.$$
In particular, $\operatorname{L}_\mathbf{z}$ attains its norm at $(x_1,\ldots, x_d)$ and $\langle \mathbf{x}, \mathbf{z}\rangle =|\lambda|^2=\Vert \operatorname{L}_\mathbf{z}\Vert^2$.

In \cite{carando2019symmetric}, the author together with Daniel Carando proved the following result, along with Lemma \ref{best rank 1 equivalence}  will be the main tool for some of the results of this section.
\begin{thm}[Carando, R.]\label{teo with carando}  Let $d\ge 3$ and $\mathbf{x}_1,\ldots,\mathbf{x}_d$ be norm one vectors on a  Hilbert space $\mathcal H$ over $\zK$.
There exists a symmetric $k$-linear form $T$ on $\mathcal H$  attaining its norm at $(\mathbf{x}_1,\ldots, \mathbf{x}_k)$ if and only if $$\dim (\operatorname{span}\{\mathbf{x}_1,\ldots,\mathbf{x}_d\})= \begin{cases} 1 \quad\quad\quad \text{ if }\zK=\zC \\ 1 \text{ or }2 \quad\text{ if }\zK=\zR.\end{cases}
$$
\end{thm}

The equivalence given in Lemma \ref{best rank 1 equivalence} and Banach's result \cite{banach1938homogene} implies that the best rank-$1$ approximation of a symmetric tensor can always be chosen to be symmetric. For the real case, in \cite{friedland2013best, friedland2014number} the authors proved that for almost every tensor in $\bigotimes^{d,s} \zR^n$ any best rank-1 approximation has to be symmetric. Moreover, in \cite[Section 4]{friedland2013best} some information about symmetric tensors that have a non-symmetric best rank-$1$ approximation is given. The next result gives some more information on the best rank-$1$ approximation of a real symmetric tensor.

\begin{prop}\label{prop best rank 1}  Let $\mathbf{z} \in \bigotimes^d \zR^n$ be a symmetric tensor and $\mathbf{x}=\lambda x_1\otimes\cdots \otimes x_d,$ with  $x_1, \ldots ,x_d\in \zR^n$ norm one vectors and $\lambda\in \zR$, a best tensor rank-$1$ approximation of $\mathbf{z}$ for the Hilbert-Schmidt norm. If we define $\mathcal{H}=\operatorname{span}\{x_1,\ldots, x_d\}$, then
$$\dim(\mathcal{H})\leq 2.$$
Moreover, if $\dim(\mathcal{H})=2$ then the multilinear form $\operatorname{L}_\mathbf{z}$ restricted to $\mathcal{H}$ is the only symmetric multilinear form on $\mathcal{H}$ with norm  $|\lambda|$ such that $\operatorname{L}_\mathbf{z}(x_1,\ldots, x_d)= \lambda$.
\end{prop}

\begin{proof}
We know that $\operatorname{L}_\mathbf{z}$ attains its norm at $(x_1,\ldots, x_d)$. Since $\mathbf{z}$ is a symmetric tensor, $\operatorname{L}_\mathbf{z}$ is a symmetric multilinear form. Then, by  Theorem \ref{teo with carando}, we conclude that 
$$\dim(\mathcal{H})\leq 2.$$

The second part of the theorem is just   \cite[Lemma 1.6]{carando2019symmetric} applied to the norm one multilinear form $\displaystyle{\frac{\operatorname{L}_\mathbf{z}}{\lambda}}$.
\end{proof}

\begin{rem}
Notice that in the previous theorem, $\dim(\mathcal{H})=1$ if and only if $\mathbf{x}$ is a symmetric tensor.
\end{rem}

For a Hilbert space $\mathcal{H}$ let $\mathcal{L}_s(^d\mathcal{H})$ be the spaces of symmetric $d$-linear forms on $\mathcal{H}$ with the usual uniform norm. As an application of Proposition \ref{prop best rank 1} we obtain a geometric characterization of the symmetric tensors with a non-symmetric best rank-$1$ approximation.

\begin{cor}\label{coro best rank 1} Let $\mathbf{z} \in \bigotimes^d \zR^n$ be a symmetric tensor. Then, for the Hilbert-Schmidt norm the tensor  $\mathbf{z}$ has a non-symmetric best tensor rank-$1$ approximation $\mathbf{x}=\lambda x_1\otimes\cdots \otimes x_d,$ with $\lambda\in \zR$ and  $x_1, \ldots ,x_d\in \zR^n$ norm one vectors, if an only if  the multilinear form $\frac{ \operatorname{L}_\mathbf{z}}{\lambda}$ restricted to $\mathcal{H}=\operatorname{span}\{x_1,\ldots, x_d\}$ is an exposed point of the unit ball of the space $\mathcal{L}_s(^d\mathcal{H})$, exposed by the linear function $f\in [\mathcal{L}_s(^d\mathcal{H})]^*$ defined as
$$f(\operatorname{L})=\operatorname{L}(x_1,\ldots, x_d).$$
\end{cor}
\begin{proof}
Assume that  $\frac{ \operatorname{L}_\mathbf{z}}{\lambda}$ is an exposed point of the unit ball. Then it has norm one. In particular 
$$\varepsilon(\mathbf{z})= \Vert \operatorname{L}_\mathbf{z}\Vert =|\lambda|.$$
The fact that
$$f\left(\frac{ \operatorname{L}_\mathbf{z}}{\lambda}\right)=1$$
implies that $\operatorname{L}_\mathbf{z}(x_1,\ldots,x_n)=\lambda$. Then, by Lemma~\ref{best rank 1 equivalence}, we have that $\lambda x_1\otimes\cdots \otimes x_d$ is a best tensor rank-$1$ approximation of $\mathbf{z}$.

On the other hand, if we assume that $\lambda x_1\otimes\cdots \otimes x_d$ is a best tensor rank-$1$ approximation of $\mathbf{z}$, by Lemma~\ref{best rank 1 equivalence}, $\operatorname{L}_\mathbf{z}(x_1,\ldots,x_n)=\lambda$, which implies that 
$$f\left(\frac{ \operatorname{L}_\mathbf{z}}{\lambda}\right)=1.$$
By Proposition \ref{prop best rank 1}, we know that there is no other symmetric multilinear form $\operatorname{L}$ such that $\operatorname{L}(x_1,\ldots,x_n)=1$. Otherwise $\lambda \operatorname{L}$ would be another symmetric multilinear form with norm  $|\lambda|$ and $\lambda\operatorname{L}(x_1,\ldots, x_d)=\lambda$. Therefore, there is no other symmetric multilinear form $\operatorname{L}$ such that 
$$f(\operatorname{L})=1.$$
That is, $f$ exposes  $\frac{ \operatorname{L}_\mathbf{z}}{\lambda}$.
\end{proof}

As the next proposition shows, if we are working on $\zR^2$ and we have a very particular non-symmetric best rank$-1$ approximation of a tensor $\mathbf{z}$, Proposition~\ref{prop best rank 1} allows to fully recover $\mathbf{z}$. But first, let us introduce some notation that will come in handy in the rest of this section.

\begin{notation}\label{notation multilinear} For a multilinear form $\operatorname{L}$, and vectors $v,w$, we write
$$\operatorname{L}(v^j, w^{d-j})= L(\underbrace{v, \cdots , v}_{j\,\,\, times}, \underbrace{w, \cdots, w}_{d-j\,\,\, times}).$$
\end{notation}

\begin{prop}\label{prop for proc} Let $v, w\in \zR^2$ be orthonormal vectors, $\mathbf{z}\in \bigotimes^d \zR^{2}$ a symmetric tensor and $1\leq j\leq d-1$. Then
$$\mathbf{x}= (\otimes^j v)\otimes (\otimes^{d-j} w):=\underbrace{v\otimes \cdots \otimes v}_{j\,\,\, times}\otimes \underbrace{w\otimes \cdots\otimes w}_{d-j\,\,\, times}$$ is a best rank-$1$ approximation of $\mathbf{z}$  if and only if
$$\mathbf{z}=
\left\{ \begin{array}{ccc}
\displaystyle{ (-1)^{\frac{j}{2}}\sum_{l=0}^{\left[\frac{d}{2}\right]} }\binom{d}{2l}(-1)^l (\vee^{2l} v)\vee (\vee^{d-2l} w) & & \text{if $j$ is even}  \\
\\
\displaystyle{(-1)^{\frac{j+1}{2}}\sum_{l=0}^{\left[\frac{d-1}{2}\right]} }\binom{d}{2l+1}(-1)^l (\vee^{2l+1} v)\vee (\vee^{d-2l-1} w)&  & \text{if $j$ is odd.} 
\end{array}\right.
$$
Where $(\vee^k v)\vee (\vee^{d-k} w):=\underbrace{v \vee \cdots \vee v}_{k\,\,\, times}\vee \underbrace{w\vee \cdots\vee w}_{d-k\,\,\, times}$ and $[\cdot]$ is the integer part.
\end{prop}
\begin{proof}
Suppose that $j$ is even and that
$$\mathbf{z}=(-1)^{\frac{j}{2}}\sum_{l=0}^{\left[\frac{d}{2}\right]} \binom{d}{2l}(-1)^l (\vee^{2l} v)\vee (\vee^{d-2l} w).$$
To see that $\mathbf{x}$  is a best rank-$1$ approximation of $\mathbf{z}$, by Lemma \ref{best rank 1 equivalence}, we need to prove that $\varepsilon(\mathbf{z})=1$ and that $\langle \mathbf{x}, \mathbf{z}\rangle =1.$ For the first part, notice that if we identify $v, w$ with the canonical basis, the polynomial $P_\mathbf{z}$ is, up to sign, essentially one of the polynomials from \cite[Lemma 4.2]{carando2019symmetric}. Both of which have norm one, therefore 
$$\varepsilon(\mathbf{z})=\Vert P_\mathbf{z}\Vert =1.$$ The second part is a straightforward computation. 
\begin{eqnarray*}
\langle  \mathbf{x}, \mathbf{z}\rangle &=&  (-1)^{\frac{j}{2}}\sum_{l=0}^{\left[\frac{d}{2}\right]} \binom{d}{2l}(-1)^l \left\langle  (\otimes^j v)\otimes (\otimes^{d-j} w), (\vee^{2l} v)\vee (\vee^{d-2l} w) \right\rangle \\
&=& (-1)^{\frac{j}{2}} \binom{d}{j} (-1)^\frac{j}{2} \left\langle (\otimes^j v)\otimes (\otimes^{d-j} w), (\vee^{j} v)\vee (\vee^{d-j} w)\right\rangle=1. \
\end{eqnarray*}
The case $j$ odd is similar.

On the other hand, if we assume that  $\mathbf{x}$  is a best rank-$1$ approximation of $\mathbf{z}$, by Lemma \ref{best rank 1 equivalence} and Proposition \ref{prop best rank 1}, we have that  $\operatorname{L}_\mathbf{z}$ is the only norm one symmetric multilinear such that 
$$\operatorname{L}_\mathbf{z}(v^j, w^{d-j})=1.$$
Suppose that $j$ is even. By the computations done in the first part we know that the symmetric multilinear form $\operatorname{L}$ associated to  the tensor
$$(-1)^{\frac{j}{2}}\sum_{l=0}^{\left[\frac{d}{2}\right]} \binom{d}{2l}(-1)^l (\vee^{2l} v)\vee (\vee^{d-2l} w)$$
is a norm one symmetric multilinear such that  
$$\operatorname{L}(v^j, w^{d-j})=1.$$
Therefore $\operatorname{L}=\operatorname{L}_\mathbf{z}$ and
$$\mathbf{z}=(-1)^{\frac{j}{2}}\sum_{l=0}^{\left[\frac{d}{2}\right]} \binom{d}{2l}(-1)^l (\vee^{2l} v)\vee (\vee^{d-2l} w).$$
The case $j$ odd is analogous.
\end{proof}

Notice that in the previous proposition, each term in the decomposition
$$\mathbf{u}=\sum_{l=0}^{\left[\frac{d}{2}\right]} \binom{d}{2l}(-1)^l (\vee^{2l} v)\vee (\vee^{d-2l} w)$$
corresponds to a best rank-$1$ approximation of $\mathbf{u}$. For each $l=0,\ldots, \left[\frac{d}{2}\right]$, the tensor
$$(-1)^l(\otimes^{2l} v)\otimes (\otimes^{d-2l} w)$$
is a  best rank-$1$ approximation of $\mathbf{u}$, and it is non-symmetric if $2l\neq 0, d$.

Next, we focus on the problem of recovering as much information as possible on a symmetric tensor $\mathbf{z}$ using a non-symmetric best rank$-1$ approximation. That is, we are interested in results similar to Proposition \ref{prop for proc}. First, let us show that we can never fully recover the tensor $\mathbf{z}$ if we are on a dimension greater than two.

\begin{prop} Let $\mathcal{H}$ be a Hilbert space of dimension at least $3$ and $d>1$. Then, for any set of norm one vectors $x_1,\ldots,x_d$ either no symmetric $d$-linear form of norm one attains its norm at them or infinitely many do.
\end{prop}

\begin{proof}
Let us assume some symmetric $d$-linear form of norm one $\operatorname{L}$ attains its norm at $x_1,\ldots,x_d$. By Theorem \ref{teo with carando} we have that 
$$\mathcal{H}_1=\operatorname{span}\{x_1,\ldots, x_d\}$$
is a subspace of dimension at most 2. Consider $\operatorname{P}:\mathcal{H} \rightarrow \mathcal{H}_1$ the orthogonal projection and $w$ a norm one vector orthogonal to  $\mathcal{H}_1$. Such $w$ exists because the dimension of $\mathcal{H}$  is at least 3. Then, for any $-1\leq a\leq 1$ we claim that
$$\operatorname{L}_a(y_1,\ldots, y_d):= \operatorname{L}(\operatorname{P}(y_1),\ldots,\operatorname{P}(y_d)) +a\, \langle y_1, w \rangle\cdots \langle y_d, w \rangle$$ is a norm one symmetric $d$-linear form that attains its norm at $x_1,\ldots,x_d$. Since
$$\operatorname{L}_a(x_1,\ldots,x_d) =\operatorname{L}(x_1,\ldots,x_d),$$
we only need to show that $\operatorname{L}_a$ has norm one. For any norm one vector $y$, we have that
\begin{eqnarray*}
|\operatorname{L}_a(y,\ldots,y)| &=& |\operatorname{L}(\operatorname{P}(y),\ldots,\operatorname{P}(y)) +a\, \langle y, w \rangle^d| \\
&\leq & \Vert P(y)\Vert^d + |\langle y, w \rangle|^d \\
&\leq & (\Vert P(y)\Vert^2 + |\langle y, w \rangle|^2)^\frac{d}{2}\\
&= & \Vert y \Vert^d=1.\
\end{eqnarray*}
Since $\operatorname{L}_a$ is symmetric, using Banach's result \cite{banach1938homogene}, this gives $\Vert \operatorname{L}_a \Vert = 1$.
\end{proof}

 In terms of tensors, this lemma states that in dimension greater than two, a tensor $x_1\otimes \ldots \otimes x_d$  is not the best rank$-1$ approximation of any symmetric tensor, or is the best  rank$-1$ approximation  of infinitely many symmetric tensors. Therefore, in this context, we can never expect to fully recover a symmetric tensor $\mathbf{z}$ using only a best rank$-1$ approximation. This remark is important because of Theorem \ref{teo recover} below.

Our main goal in this part of the paper is to show the following result. This result allows to partially recover the values of multilinear function associated with a tensor $\mathbf{z}$ using only its non-symmetric best  tensor rank$-1$ approximation.

\begin{thm}\label{teo recover} Given a symmetric tensor $\mathbf{z} \in \bigotimes^{d} \zR^n,$ and a non-symmetric best rank$-1$ approximation $x_1\otimes \cdots \otimes x_d$ we can recover the values of the multilinear form $\operatorname{L}_\mathbf{z}$ restricted to the space $\mathcal{H}=\operatorname{span}\{x_1,\ldots, x_d\}$. In particular, if $n=2$ we can fully recover $\mathbf{z}$ using only the vectors $x_1\otimes \cdots \otimes x_d$
\end{thm}

To prove this theorem, we will need several auxiliary results. Below, Lemmas \ref{lem aux0} and \ref{lem aux1} are essentially a consequence of the following fact. If $x,y\in\zR^2$ are norm one vectors which are linearly independent and $B:\zR^2\times \zR^2\rightarrow \zR$ is a norm one bilinear form such that $B(x,y)=1$, then the eigenvalues (of the corresponding matrix) are $1$ and $-1$. Therefore, if $\{x_1, x_2\}$ is an orthonormal basis of eigenvectors corresponding to the eigenvalues $1$ and $-1$ respectively, it follows that we must have
$$x=\alpha x_1+\beta x_2\,\,\,\text{ and }\,\,\, y=\alpha x_1-\beta x_2,$$
for some $\alpha, \beta\in \zR$. Fore more details on this argument see \cite[Lemma 2.1]{carando2019symmetric}.

\begin{lem}\label{lem aux0} Let $L:\zR^n \times \zR^n\rightarrow \zR$ be a norm one symmetric bilinear form and $x,y\in \zR$ be linearly independent norm one vectors such that
$$\operatorname{L}(x,y)=1.$$
Then $f_1=\frac{x+y}{\Vert x+y \Vert}$, $f_2=\frac{x-y}{\Vert x-y \Vert}$ are orthonormal vectors such that
\begin{align*} 
\operatorname{L}(f_1,f_1)&=1\\
\operatorname{L}(f_2,f_2)&=-1.
\end{align*}
\end{lem}

\begin{proof}
The orthogonality can be easily checked with an straight forward computation
\begin{eqnarray*}
\langle x+y,x-y\rangle &=& \Vert x\Vert^2+ \langle y,x\rangle -\langle x,y\rangle -\Vert y\Vert^2\\
&=& 1^2-1^2=0.\
\end{eqnarray*}

Since $\operatorname{L}(x,y)=1$, we have
\begin{eqnarray*}
1 &=&  \operatorname{L}(x,y) \\
&=& \frac{1}{4}\left( \operatorname{L}(x+y,x+y) -\operatorname{L}(x-y,x-y)\right)\\
&\leq & \frac{1}{4}\left( \Vert x+y\Vert^2 + \Vert x-y\Vert^2 \right)\\
&=& \frac{1}{4}\left( 2\Vert x\Vert^2 + 2\Vert y\Vert^2 \right)=1.\
\end{eqnarray*}
Thus
\begin{align*}
\operatorname{L}(x+y,x+y) &= \Vert x+y\Vert^2\\
\operatorname{L}(x-y,x-y)&=- \Vert x-y\Vert^2,
\end{align*}
which concludes the proof.
\end{proof}

In what follows, if $\alpha<0$, rotating a vector $x$ by an angle $\alpha$ in one direction means to rotate it $|\alpha|$ in the opposite direction.

\begin{lem}\label{lem aux1} Let $\operatorname{L}:\zR^n \times \zR^n\rightarrow \zR$ be a norm one symmetric biliniear form and $x,y\in \zR$ be linearly independent norm one vectors such that
$$\operatorname{L}(x,y)=1.$$
Take any $\alpha\in \zR$. If $\tilde  x$ is the vector obtained by rotating $x$  an angle $\alpha$ on the plane $\operatorname{span}\{x,y\}$ in any direction, and $\tilde  y$ is the vector obtained by rotating $y$  by  angle $\alpha$ in the opposite direction, then
$$\operatorname{L}(\tilde x,\tilde y)=1.$$
\end{lem}

\begin{proof}
Since we are working with rotations on the plane $\operatorname{span}\{x,y\}$, there is no harm in assuming we are on $\zR^2$. By \cite[Lemma 2.1]{carando2019symmetric}, there is an orthonormal basis $F=\{\mathbf{f}_1,\mathbf{f}_2\}$ of $\zR^2$, such that
$$[L]_F=
\left(\begin{array}{lr}
1 	& 	 0	 \\
0	&	-1		\\
\end{array}\right).$$

If $\theta$ is the angle between $x$ and $f_1$ then  (replacing $\theta$ by $-\theta$ if needed) we have
$$[x]_F= (\cos (\theta),\sin (\theta)).$$ 
On the other hand, the condition $\operatorname{L}(x,y)=1$ gives
$$[y]_F= (\cos (\theta),-\sin (\theta))=(\cos (-\theta),\sin (-\theta)).$$
Therefore, we either have
$$[\tilde x ]_F= (\cos (\theta + \alpha),\sin (\theta+\alpha))$$
$$[\tilde y ]_F= (\cos (-\theta - \alpha),\sin (-\theta-\alpha))=(\cos (\theta + \alpha),-\sin (\theta+\alpha))$$
or
$$[\tilde x ]_F= (\cos (\theta - \alpha),\sin (\theta-\alpha))$$
$$[\tilde y ]_F= (\cos (-\theta + \alpha),\sin (-\theta+\alpha))=(\cos (\theta - \alpha),-\sin (\theta-\alpha)).$$
In both cases $\operatorname{L}(\tilde x,\tilde y)=1$.
\end{proof}

In what follows, we use again Notation \ref{notation multilinear}.

\begin{lem}\label{lem aux2} For $k,l \in \zN$, let $\operatorname{L}:\zR^n \times \cdots \times \zR^n\rightarrow \zR$ be a norm one symmetric $(k+l)$-linear form and $x,y\in \zR$ be linearly independent norm one vectors with
$$\operatorname{L}(x^k,y^l)=1.$$
Take any $\alpha$ such that
$$\frac{\theta-\pi}{k+l} \leq \alpha \leq  \frac{\theta}{k+l},$$
where $\theta$ is the angle between $x$ and $y$. If $\tilde  x$ is the vector obtained by rotating $x$ an angle $l \alpha$ on the plane $\operatorname{span}\{x,y\}$ in the direction of $y$, and $\tilde  y$ is the vector obtained by rotating $y$ an angle $k \alpha$ in the opposite direction, then
$$\operatorname{L}(\tilde x^k,\tilde y^l)=1.$$
\end{lem}

\begin{proof} We only give a sketch of the proof to avoid some cumbersome notation and to be as simple and clear as possible. We have $k$ copies of $x$ and $l$ copies of $y$. Our objective is to see that we can rotate each copy of $x$ an angle of $l \alpha$ and rotate each copy of $y$ an angle of $k \alpha$ in the opposite direction.

Take one copy of $x$ and $y$. By fixing all the other variables, Lemma \ref{lem aux1} tells us that we can rotate $x$ an angle $\alpha$ and $y$ an angle $\alpha$ in the opposite direction. Repeat this for the same $x$ and a different $y$. Now we have rotated that copy of $x$ an angle $2\alpha$. Repeat this procedure with the same $x$ for the rest of the copies of $y$. By the end, we have rotated this first copy of $x$ an angle $l \alpha$ and each $y$ has been rotated an angle $\alpha$ in the opposite direction.

Next, do the same steps for another copy of $x$. After that, the second copy of $x$ has been rotated at an angle $l \alpha$ and each $y$ has been rotated an angle $\alpha$ in the opposite direction for the second time. That is, each $y$ has been rotated at an angle $2\alpha$ in the opposite direction. Repeating this for all the copies of $x$ gives the desired result.

The condition $$\frac{\theta-\pi}{k+l} \leq \alpha \leq  \frac{\theta}{k+l},$$
is there to make sure that every time we apply Lemma \ref{lem aux1} we are working with linearly independent vectors. Indeed, in each step we work with a copy of $x$ that has been rotated $i$ times and a copy of $y$ that has been rotated $j$ times, with $i<l$ and $j<k$. Therefore, the angle between those copies is $\theta - \alpha (i+j)$, and we have
$$0< \theta - \frac{\theta}{k+l}  (i+j) \leq \theta - \alpha (i+j) \leq \theta - \frac{\theta-\pi}{k+l}  (i+j)< \pi.$$
\end{proof}

\begin{rem} Notice that in the previous lemma the condition 
$$\frac{\theta-\pi}{k+l} \leq \alpha \leq  \frac{\theta}{k+l}$$
may be replaced for any condition that assure that $x$ rotated an angle $i\alpha$ is linearly independent with $y$ rotated an angle $j\alpha$ in the opposite direction. For example, one could take $\alpha$ such that for any $n\in \{0,\ldots, k+l-2\}$ the number $\theta - \alpha n$ is not  congruent to 0 modulo $\pi$. Which is a more general --but less clear-- condition.
\end{rem}

Finally, we are in conditions to prove Theorem \ref{teo recover}.

\begin{proof}[Proof of Theorem \ref{teo recover}]  To simplify the proof, we may assume $x_1,\ldots, x_d$ are norm one vectors. In particular, by Lemma \ref{best rank 1 equivalence}, this  implies that $\Vert \operatorname{L}_\mathbf{z}\Vert=1$. If this is not the case, we simply divide each $x_i$ by its norm, and this procedure will give $\operatorname{L}_\mathbf{z}$ divided by its norm.

In order to prove this result, we give a three-step algorithm to obtain orthonormal vectors $v,w\in \operatorname{span}\{x_1,\ldots, x_d\}$ such that $\operatorname{L}_\mathbf{z}(v^2,w^{d-2})=\pm 1$. This, combined with Proposition \ref{prop for proc}, shows that $\operatorname{L}_\mathbf{z}$ on $\operatorname{span}\{x_1,\ldots, x_d\}$ is the multilinear form associated to the tensor
\begin{equation}\label{eq associated}
 -\sum_{l=0}^{\left[\frac{d}{2}\right]} \binom{d}{2l}(-1)^l (\vee^{2l} v)\vee (\vee^{d-2l} w),
 \end{equation}
if $\operatorname{L}_\mathbf{z}(v^2,w^{d-2})=1$. In the case $\operatorname{L}_\mathbf{z}(v^2,w^{d-2})=-1$, remove the minus sign from the formula \eqref{eq associated}.

\textbf{Step I:} Reorder $x_1,\ldots, x_d$ to have that $x_1$ and $x_2$ are linearly independent.

\textbf{Step II:} The second step is based on a constructive proof of Banach's result given in \cite[Theorem 2.1]{pappas2007norm}. Apply Lemma \ref{lem aux2} to the bilinear form $\operatorname{L}_\mathbf{z}(x_1,\ldots, x_{d-2},\, \cdot \,, \, \cdot \,)$ obtained from fixing the first $d-2$ variables, and to the vectors $x_{d-1}, x_d$. By picking $\alpha=\frac{\theta}{2}$, where $\theta$ is the angle between the vectors, we obtain $y_1$ such that 
$$\operatorname{L}(x_1,\ldots, x_{d-2},y_1^2)=1.$$
If $x_{d-1}, x_d$ are linearly dependent just take $y_1=x_d$. In this case, if $x_d=-x_{d-1}$, also replace $x_1$ by $-x_1$ to compensate the sign.

Next do the same to the trilinear form $\operatorname{L}(x_1,\ldots, x_{d-3}, \, \cdot \, ,\, \, \cdot \, , \, \cdot \,)$ and the vectors $x_{d-2}$ and $y_1$. If we again call $\theta$ the angle between the vectors, this time choose $\alpha=\frac{\theta}{3}$. 

Continuing this procedure we find a norm one vector, let us call it $y$, such that

$$\operatorname{L}_\mathbf{z}(x_1,x_2,y^{d-2})=1.$$

\textbf{Step III:} By fixing variables, $\operatorname{L}(\, \cdot \, , \, \cdot \, , y^{d-2})$ is a bilinear symmetric form attaining its norm at $(x_1,x_2)$. Combining \textbf{Step I} and  Lemma \ref{lem aux0}, we know that $f_1=\frac{x_1+x_2}{\Vert x_1+x_2 \Vert}$, $f_2=\frac{x_1-x_2}{\Vert x_1-x_2 \Vert}$ are orthonormal vectors such that
\begin{align*} 
\operatorname{L}(f_1^2,y^{d-2})&=1\\
\operatorname{L}(f_2^2,y^{d-2})&=-1.
\end{align*}

If $f_1$ is linearly dependent with $y$, then $y$ is orthonormal to $f_2$. In this case we take $v=f_2$, $w=y$, and  have
$$\operatorname{L}(v^2,w^{d-2})=- 1.$$

If $f_1$ is linearly independent with $y$, we may apply Lemma \ref{lem aux2} to these vectors to obtain vectors $\tilde{f}_1$ and $\tilde{y}$ such that 
$$\operatorname{L}(\tilde{f}_1^2,\tilde{y}^{d-2} )=1.$$
Choosing $\alpha$ such that
$$\theta - d \alpha =\frac{\pi}{4}.$$
will ensure that  $\tilde{f}_1$ and $\tilde{y}$ are orthonormal. Then we take $v=\tilde{f}_1$ and $w=\tilde{y}$.
\end{proof}

Now let us turn our attention to the field of complex numbers $\zC$. For complex tensors, the situation is quite different from the one studied above. With a similar proof than the one from Proposition \ref{prop best rank 1}, we have the following.

\begin{prop}\label{prop best rank 1 comp}  For $d>2$, let $\mathbf{z} \in \bigotimes^d \zC^n$ be a symmetric tensor and $\mathbf{x}=x_1\otimes \cdots \otimes x_d$ a best tensor rank-$1$ approximation of $\mathbf{z}$ for the Hilbert-Schmidt norm. If we define $\mathcal{H}=\operatorname{span}\{x_1,\ldots, x_d\}$, then
$$\dim(\mathcal{H}) = 1.$$
Equivalently, $\mathbf{x}$ is symmetric.
\end{prop}

In the real case, for almost every symmetric tensor a best rank-$1$ approximation has to be symmetric (see \cite[Theorem 2]{friedland2014number}). In comparison, Proposition \ref{prop best rank 1 comp}  establishes that in the complex case for every symmetric tensor a best tensor rank-$1$ approximation has to be symmetric, as long as $d>2$. To check that this condition  is necessary, consider $\mathbf{z}\in\zC^2 \bigotimes \zC^2$ given by
$$\textbf{z}=e_1\otimes e_1 + e_2\otimes e_2,$$
where $e_1, e_2\in \zC^2$ is the canonical basis. By Lemma \ref{best rank 1 equivalence}, the tensor $$\mathbf{x}= \left(\tfrac{1}{\sqrt{2}},\tfrac{i}{\sqrt{2}}\right)\otimes \left(\tfrac{1}{\sqrt{2}},\tfrac{-i}{\sqrt{2}}\right)$$ is a best tensor rank-$1$ approximation of $\mathbf{z}$, which is not symmetric.

From a practical point of view this implies that any algorithm designed to find a best rank$-1$ approximation of a tensor $\mathbf{z}$, can be restricted only to symmetric tensors without excluding any possible solution (provided that $\mathbf{z}$ itself is symmetric).

\section{Nuclear rank}\label{sec nuclear rank}

Recall that the \emph{nuclear rank}, introduced in \cite{friedland2014nuclear}, for  a tensor $\mathbf{z} \in \bigotimes^d \zK^{n_i}$ is defined as
$$\operatorname{rank}_\pi(\mathbf{z}):= \min\left\{r\in \zN: \mathbf{z}=\sum_{i=1}^r z_1^i\otimes \cdots \otimes z_d^i, \pi(\mathbf{z})=\sum_{i=1}^r \Vert z_1^i \Vert \cdots \Vert z_d^i\Vert\right\}.$$
For any tensor in $\bigotimes^d \zK^{n_i}$ there is always a representation that gives its nuclear norm (see for example \cite[Proposition 3.1]{friedland2014nuclear}), thus the nuclear rank is well defined.

In this section, we study the representations of symmetric tensors that give the nuclear norm and are used in the nuclear rank definition. That is, for a symmetric tensor $\mathbf{z}$, we investigate the representations $\mathbf{z}=\sum_{i=1}^r z_1^i\otimes \cdots \otimes z_d^i$ such that
$$\pi(\mathbf{z})=\sum_{i=1}^r \Vert z_1^i \Vert \cdots \Vert z_d^i\Vert.$$
Our objective is to prove the following result.
\begin{thm}\label{teo rep for sym} Let $\mathbf{z} \in \bigotimes^d \zK^n$ be a symmetric tensor  with $d>2$. If
$$\mathbf{z}=\sum_{i=1}^r z_1^i\otimes \cdots \otimes z_d^i\,\,\,\text{ with }\,\,\, \pi(\mathbf{z})=\sum_{i=1}^r \Vert z_1^i \Vert \cdots \Vert z_d^i\Vert,$$
then for each $i=1,\ldots, n$ we have 
$$\dim (\operatorname{span}\{z_1^i,\ldots,z_d^i\})= \begin{cases} 1 \quad\quad\quad \text{ if }\,\,\zK=\zC \\ 1 \text{ or }2 \quad\text{ if }\,\,\zK=\zR.\end{cases}
$$
\end{thm}

Finding representations of a tensor that gives its projective norm can be a very difficult task. The previous theorem limits considerable the space of possible representations. In particular, in the complex case, in order to compute the projective norm and the projective rank of a tensor $\mathbf{z}$ we can restrict ourselves to representations of the form
$$\mathbf{z}=\sum_{i=1}^r z^i\otimes \cdots \otimes z^i.$$

To prove this, we need  a result  analogous  to \cite[Lemma 4.1]{friedland2014nuclear} for symmetric tensors.

\begin{lem}\label{lem nuc rank} Let $\mathbf{z} =\sum_{i=1}^r z_1^i\otimes \cdots \otimes z_d^i\in \bigotimes^d \zK^n$ be a symmetric tensor. Then $$\pi(\mathbf{z})=\sum_{i=1}^r \Vert z_1^i \Vert \cdots \Vert z_d^i\Vert,$$
if and only if there is a norm one symmetric multilinear form $\operatorname L$ such that 
$$\operatorname L(z_1^i,\ldots, z_d^i)=\Vert z_1^i \Vert \cdots \Vert z_d^i\Vert\,\,\,\text{ for }\,\,\, i=1,\ldots, r.$$
\end{lem}
\begin{proof} Assume that $$\pi(\mathbf{z})=\sum_{i=1}^r \Vert z_1^i \Vert \cdots \Vert z_d^i\Vert.$$
and let us find the symmetric multilinear form. Since $\bigotimes^d \zK^n$ is finite dimensional, there is a tensor $\mathbf{u}\in \bigotimes^d \zK^n$  such that $\varepsilon(\mathbf{u})=1$ and
$$|\langle \mathbf{z}, \mathbf{u}\rangle|= \pi(\mathbf{z})\varepsilon(\mathbf{u})=\pi(\mathbf{z}).$$
First we are going to show that we can take this a symmetric tensor, and that we can also assume $|\langle \mathbf{z}, \mathbf{u}\rangle|=\langle \mathbf{z}, \mathbf{u}\rangle$. If $\mathbf{u}$ is not symmetric, consider the symmetric tensor $\mathbf{w}:=\sigma(\mathbf{u})$. Using that $\sigma: (\bigotimes^d \zK^n, \varepsilon)\rightarrow  (\bigotimes^{d,s} \zK^n, \varepsilon_s)$ has norm one (see \cite[Section 3.1]{floret1997natural}), we have 
$$\varepsilon(\mathbf{w})=\varepsilon_s(\mathbf{w})=\varepsilon_s(\sigma(\mathbf{u})) \leq \varepsilon(\mathbf{u}) =1.$$
Since $\mathbf{z}$ is symmetric, for any permutation $\eta\in S_d$, and any set of vectors $x_1,\ldots, x_d\in \zK^n$, we have
$$\langle \mathbf{z}, x_1\otimes \cdots \otimes x_d \rangle =\langle \mathbf{z}, x_{\eta(1)}\otimes \cdots \otimes x_{\eta(d)} \rangle.$$
In particular
$$ \langle \mathbf{z}, \mathbf{u} \rangle=\langle \mathbf{z}, \sigma(\mathbf{u}) \rangle=\langle \mathbf{z}, \mathbf{w} \rangle.$$
Then we have
$$    \pi(\mathbf{z})=|\langle \mathbf{z}, \mathbf{u} \rangle|= |\langle \mathbf{z}, \mathbf{w} \rangle| \leq  \pi(\mathbf{z}) \varepsilon(\mathbf{w})  \leq \pi(\mathbf{z}).$$
Hence, replacing $\mathbf{u}$ by $\mathbf{w}$ if needed, we may assume $\mathbf{u}$ is symmetric. By multiplying $\mathbf{u}$ with the modulus one scalar  $\frac{\langle \mathbf{z}, \mathbf{u} \rangle}{|\langle \mathbf{z}, \mathbf{u} \rangle|}$ we may also assume that $\langle \mathbf{z}, \mathbf{u} \rangle\in \zR_{>0}$. 

We claim that $\operatorname{L}_\mathbf{u}$ is the multilinear form we wanted to find. Given that $\mathbf{u}$ is symmetric, so is $\operatorname{L}_\mathbf{u}$. Since $\varepsilon(\mathbf{u})=1$, $\operatorname{L}_\mathbf{u}$ has norm one. Finally, the equation
$$\sum_{i=1}^r \Vert z_1^i \Vert \cdots \Vert z_d^i\Vert=\pi(\mathbf{z}) = \langle \mathbf{z}, \mathbf{u} \rangle = \sum_{i=1}^r \langle z_1^i\otimes \cdots \otimes z_d^i, \mathbf{u} \rangle = \sum_{i=1}^r \operatorname{L}_\mathbf{u}(z_1^i,\ldots, z_d^i)$$
and the fact that $\operatorname{L}_\mathbf{u}$ has norm one implies that 
$$\operatorname{L}_\mathbf{u}(z_1^i,\ldots, z_d^i) =  \Vert z_1^i \Vert \cdots \Vert z_d^i\Vert$$
for each $i=1,\ldots, r$.

Now let us assume that  there exist a multilinear form $\operatorname L$ as in the statement. Let $\mathbf{u}$ be the tensor associated to $\operatorname L$. That is, $\operatorname L=\operatorname{L}_\mathbf{u}.$ Then
\begin{eqnarray*}
\pi(\mathbf{z})&\leq& \sum_{i=1}^r \Vert z_1^i \Vert \cdots \Vert z_d^i\Vert = \sum_{i=1}^r \operatorname{L}_\mathbf{u}(z_1^i,\ldots, z_d^i) \\
&=&\left|\sum_{i=1}^r \operatorname{L}_\mathbf{u}(z_1^i,\ldots, z_d^i) \right|= |\langle \mathbf{z}, \mathbf{u} \rangle| \leq \pi(\mathbf{z}) \varepsilon(\mathbf{u})=\pi(\mathbf{z}).\
\end{eqnarray*}
Therefore, all the inequalities are in fact equalities, as we wanted to see.
\end{proof}

Having proved this lemma, Theorem \ref{teo rep for sym} follows from combining it with Theorem \ref{teo with carando}.

As a final remark of this section let us show that $\dim (\operatorname{span}\{z_1^i,\ldots,z_d^i\})=2$ can occur on the real case. Moreover, in the following we exhibit example in which a representation giving the nuclear norm fulfils that $\operatorname{span}\{z_1^i,\ldots,z_d^i\})$ has dimension two for every $i$.
 
\begin{exa} Let $x_1,\ldots, x_d$ be vectors on $\zR^n$ such that $\operatorname{span}\{x_1,\ldots,x_d\})$ has dimension two. Then 
$$x_1\vee\cdots \vee x_d = \frac{1}{d!}\sum_{\eta \in S_d} x_{\eta(1)}\otimes \cdots \otimes x_{\eta(d)}$$
is a symmetric tensor and
$$\pi(x_1\vee\cdots \vee x_d )= \frac{1}{d!}\sum_{\eta \in S_d} \Vert x_{\eta(1)}\Vert\cdots \Vert x_{\eta(d)}\Vert.$$

Indeed, we have
$$\Vert x_1 \Vert \cdots \Vert x_d\Vert = \pi(x_1\vee\cdots \vee x_d )\leq \frac{1}{d!}\sum_{\eta \in S_d} \Vert x_{\eta(1)}\Vert\cdots \Vert x_{\eta(d)}\Vert = \Vert x_1 \Vert \cdots \Vert x_d\Vert,$$
where the first equality follows from \cite[Proposition 1.3]{carando2019symmetric}.
\end{exa}

\section{Decomposable symmetric rank}\label{sect dec rank}

In this section we introduce the notion of decomposable symmetric rank. For  a symmetric tensor $\mathbf{z}~\in~\bigotimes^{d,s}~\zK^n$, its \emph{decomposable symmetric rank} is defined as
$$\operatorname{rank}_\sigma(\mathbf{z}):= \min\left\{r\in \zN: \mathbf{z}=\sum_{i=1}^r z_1^i\vee \cdots \vee z_d^i\right\}.$$
It is important to remark that this definition is only for  symmetric tensors.

One advantage over the tensor rank and the nuclear rank is that the decomposable symmetric rank never exceeds them.

\begin{prop}\label{prop rank} Let $\mathbf{z} \in \bigotimes^{d,s} \zK^n$ be a symmetric tensor. Then
$$\operatorname{rank}_\sigma(\mathbf{z}) \leq \operatorname{rank}(\mathbf{z})\leq \operatorname{rank}_\pi(\mathbf{z}).$$
\end{prop}
\begin{proof}
The inequality
$$\operatorname{rank}(\mathbf{z})\leq \operatorname{rank}_\pi(\mathbf{z})$$
is immediate from the definitions of tensor rank and nuclear rank. The other inequality follows from the fact that if $\mathbf{z}$ is a symmetric tensor  with a representation 
$\mathbf{z} =\sum_{i=1}^r z_1^i\otimes \cdots \otimes z_d^i,$
then
$$\mathbf{z} =\sigma(\mathbf{z}) =\sum_{i=1}^r z_1^i\vee \cdots \vee z_d^i.$$
\end{proof}

It is not hard to see that the decomposable symmetric rank can  differ from the tensor rank. Take two linear independent vectors $v,w\in \zK^n$. If we consider $\mathbf{z}=v\vee w\in \bigotimes^2\zK^n$, then
$$\operatorname{rank}_\sigma(\mathbf{z})=1 \,\,\,\text{ and } \,\,\, \operatorname{rank}(\mathbf{z})=2.$$

In the following proposition, we show that in the two-dimensional complex space $\zC^2$ every symmetric tensor has decomposable symmetric rank 1.

\begin{prop}Let $\mathbf{z} \in \bigotimes^{d} \zC^2$ be a symmetric tensor. Then
$$\operatorname{rank}_\sigma(\mathbf{z})=1.$$
\end{prop}

\begin{proof} We need to prove that there are vectors $z_1,\ldots,z_d\in \zC^2$ such that
$$\mathbf{z}=z_1\vee\cdots \vee z_d.$$
In terms of polynomials, this corresponds to find $z_1,\ldots,z_d\in \zC^2$ with
$$\operatorname{P}_\mathbf{z}(y)=\langle y, z_1\rangle \cdots \langle y, z_d\rangle.$$
Therefore, this result is equivalent to Proposition \ref{prop pol on C2} below.
\end{proof}

The following is a known result, we add a proof here for completeness.

\begin{prop}\label{prop pol on C2} Let $\operatorname P:\zC^2\rightarrow \zC$ be a $d$-homogeneous polynomial. Then, there are linear functions $\varphi_1,\ldots,\varphi_d:\zC^2\rightarrow \zC$ such that
$$\operatorname P=\prod_{i=1}^d \varphi_i.$$
\end{prop}
\begin{proof} We proceed by induction on $d$. The case $d=1$ is trivial. Let us prove the result for $d>1$, assuming it holds for $d-1$. Take $B=\{b_1,b_n\}$ a basis of $\zC^2$ such that $\operatorname P(b_1)=0$. Next we write $\operatorname P$ in function of this basis
$$\operatorname P(y)=\sum_{\substack{i_1,  i_2\in \zN_0,\\ i_1+i_2=d}} a_{i_1i_2} [y]_B^{(i_1,i_2)},$$
where $[y]_B=(t_1,t_2)$ stands for the coordinates of $y$ in the basis $B$ and $[y]_B^{(i_1,i_2)}=t_1^{i_1}t_2^{i_2}.$ Using that $\operatorname P$ is zero on $b_1$, it is easy to deduce that the coefficient $a_{d\,0}$ is zero:
$$0=\operatorname P(b_1)=\sum_{\substack{i_1, i_2\in \zN_0,\\ i_1+i_2=d}} a_{i_1i_2} (1,0)^{(i_1,i_2)}=a_{d\,0}.$$

Then we have
\begin{eqnarray*}
\operatorname P(y)&=&\sum_{\substack{i_1,i_2\in \zN_0,\\ i_1+i_2=d}}  a_{i_1i_2} [y]_B^{(i_1,i_2)}=\sum_{\substack{i_1,i_2\in \zN_0,\\ i_1+i_2=d, i_2\neq 0}}  a_{i_1i_2} [y]_B^{(i_1,i_2)}\\
&=& t_2 \sum_{\substack{i_1,i_2\in \zN_0,\\i_1+i_2=d, i_2\neq 0}} a_{i_1i_2} [y]_B^{(i_1,i_2-1)}.\
\end{eqnarray*}
Applying the inductive hypothesis to the $(d-1)$-homogeneous polynomial
$$\operatorname Q(y)=\sum_{\substack{i_1, i_2\in \zN_0,\\i_1+i_2=d, i_2\neq 0}} a_{i_1i_2} [y]_B^{(i_1,i_2-1)}$$ we obtain the desired result.
\end{proof}

\begin{rem} As said many times, two main reasons to work with low-rank tensors is that in general they are cheaper to store and computations using them require less time. As seen in Proposition \ref{prop rank}, the use of  the decomposable symmetric rank is a good alternative to reduce the capacity needed to store symmetric tensors. But as far as computations go, we have the problem that each decomposable symmetric tensor is formed by $d!$ elementary tensors. Because of this, computing 
$$ \langle z_1\vee \cdots \vee z_d, x_1\otimes\cdots\otimes x_d \rangle= \frac{1}{d!}\sum_{\eta \in S_d} \langle z_{\eta(1)}\otimes \cdots \otimes z_{\eta(d)}, x_1\otimes\cdots\otimes x_d \rangle $$
can be very costly. This implies that in the most general setting, using the decomposable symmetric rank may not be the optimal alternative to do computations. But if all the tensors involved are symmetric this is no longer a problem, since
$$\langle z_1\vee \cdots \vee z_d, x_1\vee\cdots \vee x_d\rangle = \langle z_1\otimes \cdots \otimes z_d, x_1\otimes \cdots \otimes x_d\rangle.$$
This last equality was observed on the proof of Lemma \ref{lem nuc rank}.
\end{rem}

Given a tensor $\mathbf{z}$ and   a low-rank approximation $\mathbf{x}$, when we consider the decomposable symmetric rank, we impose on $\mathbf{x}$ to be symmetric. Therefore, if $\mathbf{z}$ is far away from the space of symmetric tensors it is not possible to find a good approximation with low decomposable symmetric rank. On the opposite extreme, as the next result shows, if $\mathbf{z}$ is a symmetric tensor, a low decomposable symmetric rank approximation is at least as good as a low tensor rank approximation.

\begin{thm}\label{teo aprox dec} Let $\mathbf{z} \in \bigotimes^{d} \zK^n$ be a symmetric tensor and $\alpha(\, \cdot \,)$ be either the injective, the projective or the Hilbert-Schmidt norm. For any tensor $\mathbf{y}\in \bigotimes^{d} \zK^n$, $\mathbf{x}=\sigma(\mathbf{y})$ is a symmetric tensor with
$$\alpha(\mathbf{z}- \mathbf{x}) \leq \alpha(\mathbf{z}-\mathbf{y}).$$
Moreover, the decomposable symmetric rank of $\mathbf{x}$ is less or equal than the tensor rank and the nuclear rank of $\mathbf{y}$.
\end{thm}
\begin{proof}  To prove the first part it is enough to show that for any tensor $\mathbf{w}\in \bigotimes^{d} \zK^n$ the following inequality holds 
\begin{equation}\label{eq 1}
\alpha(\sigma(\mathbf{w}))\leq \alpha(\mathbf{w}).
\end{equation}
If we prove this, taking $\mathbf{w} = \mathbf{z}-\mathbf{x}$, we have
$$\alpha(\mathbf{z}- \mathbf{x}) \leq \alpha(\mathbf{z}-\mathbf{y}).$$

Let us prove \eqref{eq 1} for the injective norm. For this norm, by the results on \cite[Section 3.1]{floret1997natural}, we have the following
$$\varepsilon_s(\sigma(\mathbf{w}))\leq \varepsilon(\mathbf{w}).$$
As mentioned in Section \ref{sec prel}, for the particular case $\bigotimes^{d} \zK^n$, $\varepsilon_s(\, \cdot \,)$ is the restriction of the injective norm. Hence, $\varepsilon(\sigma(\mathbf{w}))=\varepsilon_s(\sigma(\mathbf{w}))$.

Next we deal with the projective norm. By the results on \cite[Section 2.3]{floret1997natural}, we have that
$$\pi_s(\sigma(\mathbf{w}))\leq c(d,\zK^n)\pi(\mathbf{w}),$$
where $c(d,\zK^n)$ is the $d-$th polarization constant of $\zK^n$ (see equation \eqref{def polarization}). For Hilbert spaces, such as $\zK^n$, the polarization constant is one. This is consequence of Banach's result  \cite{banach1938homogene} (see  also \cite[Section 2.1]{floret1997natural}). Then, \eqref{eq 1} follows from this and the fact that, since we are on a Hilbert space, $\pi(\sigma(\mathbf{w}))=\pi_s(\sigma(\mathbf{w}))$.

When considering the Hilbert-Schmidt norm, $(\bigotimes^{d} \zK^n, \operatorname{HS})$ is a Hilbert space and  $\sigma$ is the orthogonal projection onto the subspace $(\bigotimes^{d,s} \zK^n, \operatorname{HS})$. This gives \eqref{eq 1} for the Hilbert-Schmidt norm.

For the second part, notice that for any representation $\mathbf{y} =\sum_{i=1}^r y_1^i\otimes \cdots \otimes y_d^i,$
then
$$\mathbf{x} =\sigma(\mathbf{y}) =\sum_{i=1}^r y_1^i\vee \cdots \vee y_d^i.$$
Which implies that $\operatorname{rank}_\sigma(\mathbf{x}) \leq \operatorname{rank}(\mathbf{y})$. \end{proof}
As a consequence  of this result, any algorithm used to find a low tensor rank approximation $\mathbf{y}$ of a tensor $\mathbf{z}$, is suitable to find a low decomposable symmetric rank  approximation $\mathbf{x}$ as good as $\mathbf{y}$. This is true for the injective, the projective, and the Hilbert-Schmidt norms. Since $\mathbf{x}$ is given by $\sigma(\mathbf{y})$, the amount of data needed to store $\mathbf{x}$ is the same as the one needed for $\mathbf{y}$ or less. 

From the proof of Theorem \ref{teo aprox dec} we have that  $\operatorname{HS}(\mathbf{z}- \mathbf{x}) \leq \operatorname{HS}(\mathbf{z}-\mathbf{y})$ is an equality if and only if $\mathbf{z}- \mathbf{x}=\mathbf{z}-\mathbf{y}$. Therefore, if $\mathbf{y}$ is not symmetric, then $\mathbf{x}$ is a strictly better approximation than $\mathbf{y}$, in the sense that $\operatorname{HS}(\mathbf{z}- \mathbf{x}) < \operatorname{HS}(\mathbf{z}-\mathbf{y})$. In the following proposition we show that this is also true for the projective norm on $\bigotimes^{d} \zC^n$, provided that $d>2$. This proposition may have some interest on its own for the study of the projective tensor products of Hilbert spaces.

\begin{prop} Let $\mathbf{w}\in \bigotimes^{d} \zC^n$ be a non-symmetric tensor. If $d>2$ then
$$\pi(\sigma(\mathbf{w}))<\pi(\mathbf{w}).$$
\end{prop}
\begin{proof}Take a representation 
$$\mathbf{w}=\sum_{i=1}^r w_1^i\otimes \cdots \otimes w_d^i\,\,\,\text{ with }\,\,\, \pi(\mathbf{w})=\sum_{i=1}^r \Vert w_1^i \Vert \cdots \Vert w_d^i\Vert.$$
Since $\mathbf{w}$ is a non-symmetric tensor, for some $i_0$ we have that $\operatorname{span}\{w_1^{i_0}, \cdots , w_d^{i_0}\}$ has dimension greater than one. Then, by \cite[Proposition 1.3]{carando2019symmetric}, we have
$$\pi(w_1^{i_0}\vee \cdots \vee w_d^{i_0}) < \Vert w_1^{i_0}\Vert \cdots \Vert w_d^{i_0} \Vert.$$
Using that $\pi\left(w_1^i\vee \cdots \vee w_d^i \right)\leq \Vert w_1^i \Vert \cdots \Vert w_d^i\Vert$ for every $i$, and that for $i_0$ this inequality is strict we obtain the desired result:
\begin{eqnarray*}
\pi(\sigma(\mathbf{w})) &=& \pi\left( \sum_{i=1}^r w_1^i\vee \cdots \vee w_d^i \right) \\
&\leq&  \sum_{i=1}^r \pi\left(w_1^i\vee \cdots \vee w_d^i \right) \\
&<&  \sum_{i=1}^r \Vert w_1^i \Vert \cdots \Vert w_d^i\Vert \\
&=& \pi(\mathbf{w}).\
\end{eqnarray*}
\end{proof}

Below we give some examples for which this result does not hold.
\begin{itemize}
\item If $v, w\in \zC^n$ are orthonormal vectors, then $\mathbf{w}=v\otimes w$ is a non-symmetric tensor with $\pi(\sigma(\mathbf{w}))=\pi(\mathbf{w}).$
\item In $\bigotimes^d\zR^n$, by \cite[Proposition 1.3]{carando2019symmetric}, any tensor of the form $\mathbf{w}=w_1\otimes\cdots\otimes w_d$ with $$\operatorname{dim}(\operatorname{span}\{w_1,\cdots, w_d\})=2$$ is an example of a non-symmetric tensor such that $\pi(\sigma(\mathbf{w}))=\pi(\mathbf{w}).$
\item For the injective tensor norm let us build an example on $\bigotimes^2\zK^3$. Take $e_1,e_2, e_3\in \zK^3$ the canonical basis. For $t>0$ small enough, the tensor
$$\mathbf{w}=e_1\otimes e_1 + t(e_2\otimes e_3-e_3\otimes e_2)$$
has injective norm one, and $\varepsilon(\sigma(\mathbf{w}))=\varepsilon(e_1\otimes e_1)=1$.
\end{itemize}

Next, we exhibit an example showing that the decomposable symmetric rank shares one of the problems of the tensor rank: there is not always a best rank-$r$ approximation.

\begin{exa}\label{exam rank 2} On $\bigotimes^{3} \zK^6$, consider the canonical basis $\{e_1,\ldots, e_6\}$ of $\zK^6$ and the sequence of symmetric tensors
$$\mathbf{y}_n:= n\left(e_1+ \frac{1}{n} e_4\right)\vee\left(e_2+\frac{1}{n} e_5\right)\vee\left(e_3+\frac{1}{n} e_6\right)-n \,e_1\vee e_2\vee e_3.$$
Then $\mathbf{y}_n$ converges to the tensor
$$\mathbf{y}:= e_1\vee e_2\vee e_6 +e_1\vee e_3 \vee e_5+e_2\vee e_3 \vee e_4.$$
This follows from the computation
$$\mathbf{y}_n-\mathbf{y}= \frac{1}{n}([e_3\vee e_4\vee e_5] + [e_2\vee e_4\vee e_6] + [e_1\vee e_5\vee e_6]) + \frac{1}{n^2}e_4\vee e_5\vee e_6.$$
Clearly $\operatorname{rank}_\sigma(\mathbf{y}_n)\leq 2$. But, as next lemma shows, $\operatorname{rank}_\sigma(\mathbf{y})= 3.$ In particular, we have that the set
$$\{\mathbf{u}\in {\textstyle\bigotimes ^{3}} \zK^6: \mathbf{u} \text{ has decomposable rank at most } 2\}$$
is not closed. Thus, there is not always a best decomposable symmetric rank-$2$ approximation on this space. This is independent on the norm considered on $\bigotimes^{3} \zK^6$, since this is a finite dimensional space all norms are equivalent.
\end{exa}

\begin{lem} The tensor $\mathbf{y}$ from then previous example has decomposable symmetric rank $3$.
\end{lem}
\begin{proof}
We need to show that $\operatorname{rank}_\sigma(\mathbf{y})\geq 3.$ Let us assume that there are vectors $y_1,\ldots, y_6\in \zK^6$ such that
\begin{equation}\label{eq alt rep}
\mathbf{y}=y_1\vee y_2\vee y_3 + y_4\vee y_5\vee y_6
\end{equation}
and arrive to a contradiction. To do this we are going to use an auxiliary multilinear operator $E:\bigotimes^{3} \zK^6\times \zK^6 \times \zK^6\rightarrow \zK^6$ defined on elementary tensors as
$$E(x_1\otimes x_2 \otimes x_3, v, w)= \langle v, x_1\rangle \langle w, x_2\rangle x_3.$$

As an informative note, in terms of multilinear forms, this operator is just
$$E(\mathbf{u}, v, w)= \operatorname{L}_\mathbf{u}(v, w, \, \cdot \,) \in (\zK^6)^*,$$
and then identifying $\zK^6$ with its dual $(\zK^6)^*$  via the inner product.

First we are going to show that $y_1,\ldots, y_6$ is a basis of $\zK^6$. Fix a variable to obtain a bilinear operator $E(\mathbf{y}, \, \cdot \, , \, \cdot \,):\zK^6\times \zK^6\rightarrow\zK^6$. Then, by \eqref{eq alt rep}, this  bilinear operator  has its image included in $\operatorname{span}\{y_1,\ldots, y_6\}$. On the other hand, using that $\mathbf{y}:= e_1\vee e_2\vee e_6 +e_1\vee e_3 \vee e_5+e_2\vee e_3 \vee e_4$, the computations
$$E(\mathbf{y}, e_2, e_6)= \frac{1}{3!}e_1$$
$$E(\mathbf{y}, e_1, e_6)= \frac{1}{3!}e_2$$
$$E(\mathbf{y}, e_1, e_5)= \frac{1}{3!}e_3$$
$$E(\mathbf{y}, e_2, e_3)= \frac{1}{3!}e_4$$
$$E(\mathbf{y}, e_1, e_3)= \frac{1}{3!}e_5$$
$$E(\mathbf{y}, e_1, e_2)= \frac{1}{3!}e_6,$$
show that the image of $E(\mathbf{y}, \, \cdot \, , \, \cdot \,)$ is $\zK^6$. Therefore 
$\zK^6\subseteq  \operatorname{span}\{y_1,\ldots, y_6\},$ as we wanted to see.

Now take $y_1^*,\ldots, y_6^*\in \zK^6$ the dual basis of $y_1,\ldots, y_6$. That is 
$$\langle y_i^*, y_j\rangle =\delta_{ij},$$
where $\delta_{ij}$ is the Kronecker delta. At least one of the vectors of this basis is not orthogonal to $e_1$. There is no harm in assuming $y_1^*$ is not orthogonal to $e_1$. 

Consider the linear operator $E(\mathbf{y}, y_1^*, \, \cdot \,):\zK^6\rightarrow\zK^6$ obtained from fixing two variables. By equation \eqref{eq alt rep}, and the fact that
$$\langle y_1^*, y_j\rangle \neq 0 $$
if and only if $j=1$,  we have that the image of this operator is included in the two-dimensional space $\operatorname{span}\{y_2,y_3\}$. 

Since $\langle y_1^*, e_1\rangle\neq 0$, then
\begin{align*}
E(\mathbf{y}, y_1^*, e_5)&= \frac{1}{3!} \langle y_1^*, e_1\rangle e_3 +  \frac{1}{3!} \langle y_1^*, e_3\rangle e_1 \\
E(\mathbf{y}, y_1^*, e_6)&= \frac{1}{3!} \langle y_1^*, e_1\rangle e_2 +  \frac{1}{3!} \langle y_1^*, e_2\rangle e_1 \\
E(\mathbf{y}, y_1^*, e_2)&= \frac{1}{3!} \langle y_1^*, e_1\rangle e_6 +  \frac{1}{3!} \langle y_1^*, e_6\rangle e_1+  \frac{1}{3!} \langle y_1^*, e_3\rangle e_4+  \frac{1}{3!} \langle y_1^*, e_4\rangle e_3 \\
\end{align*}
implies that the image of $E(\mathbf{y}, y_1^*, \, \cdot \,)$ has dimension at least $3$, which is the desired contradiction.
\end{proof}

Although there is not always a best rank-$r$ approximation when we are considering the decomposable symmetric rank, this is not the case for $r=1$. To end this section we show that there is always a best  decomposable symmetric rank-$1$ approximation. In order to do this, we give a characterization of the best rank-$1$ approximation for the Hilbert-Schmidt norm. This will be obtained in a similar manner as it is done in Lemma \ref{best rank 1 equivalence} for the tensor rank.

\begin{lem}\label{best rank 1 equivalence dec} Let $\mathbf{z} \in \bigotimes^d \zK^{n}$ and $\mathbf{x}=\lambda x_1\vee\cdots \vee x_d,$ with $\lambda\in \zK$ and  $x_1, \ldots ,x_d$ norm one vectors, a decomposable symmetric tensor. Then $\mathbf{x}$  is a best decomposable symmetric rank-$1$ approximation of  $\mathbf{z}$ for the Hilbert-Schmidt norm if and only if 
$$\frac{|\operatorname L_\mathbf z(x_1\ldots, x_d)|}{\operatorname{HS}(x_1\vee\cdots \vee x_d)} = \max \left\{\frac{|\operatorname L_\mathbf z(u_1\ldots, u_d)|}{\operatorname{HS}(u_1\vee\cdots \vee u_d)}: \Vert u_1\Vert =\cdots =\Vert u_d\Vert=1\right\}$$
and $\lambda$ is given by the formula 
$$\lambda= \frac{\langle \mathbf{z} , x_1\vee\cdots \vee x_d\rangle}{[\operatorname{HS}( x_1\vee\cdots \vee x_d\rangle)]^2}.$$

\end{lem}
\begin{proof}
For any set of norm one vectors $u_1,\ldots, u_d$ let $$\mathbf{U}:= \operatorname{span}\{u_1\vee\cdots \vee u_d\} \subseteq \bigotimes^d \zK^{n}.$$ Consider $\operatorname P _\mathbf{U}$ the orthonormal projection on this subspace, which is given by the formula
\begin{equation}\label{ec first}
\operatorname P _\mathbf{U}(\mathbf{y})=\langle \mathbf{y} , u_1\vee\cdots \vee u_d\rangle \frac{u_1\vee\cdots \vee u_d}{[\operatorname{HS}(u_1\vee\cdots \vee u_d)]^2}.\end{equation}

Since
$$
\min  \{\operatorname{HS}\left(\mathbf{z} - t u_1\vee\cdots \vee u_d): t\in \zR\right\}=\operatorname{HS} (P _\mathbf{U^\bot}(\mathbf{z})),
$$
and the fact that
$$\operatorname{HS}(\mathbf{z})^2=[\operatorname{HS} (P _\mathbf{U}(\mathbf{z}))]^2+[\operatorname{HS} (P _\mathbf{U^\bot}(\mathbf{z}))]^2,$$
obtaining a minimum for the expression
$$\inf \{\operatorname{HS}\left(\mathbf{z} - t u_1\vee\cdots \vee u_d\right): t\in\zR, \Vert u_1\Vert =\cdots =\Vert u_d\Vert=1\}$$
is equivalent to finding a maximum for
\begin{equation}\label{caract best 1}
\sup \left\{\operatorname{HS}(\operatorname P _\mathbf{U}(\mathbf{z})):\Vert u_i\Vert =1\, \forall i\right\}=\sup \left\{\frac{|\langle \mathbf{z} , u_1\vee\cdots \vee u_d\rangle| }{\operatorname{HS}(u_1\vee\cdots \vee u_d)}: \Vert u_i\Vert =1 \,\forall i\right\}.\end{equation}

But $(u_1,\ldots, u_d)\longrightarrow \frac{|\langle \mathbf{z} , u_1\vee\cdots \vee u_d\rangle| }{\operatorname{HS}(u_1\vee\cdots \vee u_d)}$  is a continuous function, therefore it does attain its maximum over the compact set
$$\{(u_1,\ldots, u_d):\Vert u_1\Vert =\cdots =\Vert u_n\Vert=1\}.$$
This implies that the supreme on \eqref{caract best 1} is in fact maximum, and proves the  assertion regarding the vectors $x_1, \ldots ,x_d$. The statement about $\lambda$ follows from \eqref{ec first}.
\end{proof}

Although the formula on Lemma \ref{best rank 1 equivalence dec} seems rather convoluted, and may not lead to an algorithm to find a best  decomposable symmetric rank$-1$ approximation, it does ensure that such a tensor exists. Which, as pointed out in \cite{comon2008symmetric, de2008tensor}, is not a minor fact. On the other hand, as shown in Theorem \ref{teo aprox dec}, any algorithm designed to obtain a best tensor rank-$1$ approximation $\mathbf{y}$ is suitable to obtain a decomposable symmetric rank-$1$ approximation $\mathbf{x}$, which is at least as good as $\mathbf{y}$.






\begin{thebibliography}{10}




\bibitem{banach1938homogene}
S.~Banach.
\emph{{\"U}ber homogene polynome in (${L}^2$).}
Studia Mathematica {\bf 7} (1938), pp. 36--44.

\bibitem{bochnak1971polynomials} J. Bochnak and J. Siciak.
\emph{Polynomials and multilinear mappings in topological vector-spaces}.
Studia Mathematica {\bf 39} (1971), pp. 59--76.

\bibitem{carando2019symmetric}
D.~Carando and J.~T. Rodr{\'\i}guez.
\emph{Symmetric multilinear forms on {H}ilbert spaces: Where do they attain their norm?}
Linear Algebra and its Applications {\bf 563} (2019), pp. 178--192.

\bibitem{cichocki2015tensor}
A. Cichocki, D. Mandic, A-H. Phan, C. Caiafa, G.~Zhou, Q.~Zhao and L.~De~Lathauwer.
\emph{Tensor Decompositions for Signal Processing Applications: From Two-way to Multiway Component Analysis.}
IEEE signal processing magazine \textbf{32} (2015), pp. 145--163.



\bibitem{comon2008symmetric}
P. Comon, G. Golub, L. H. Lim and B. Mourrain.
\emph{Symmetric tensors and symmetric tensor rank.} 
SIAM Journal on Matrix Analysis and Applications \textbf{30} (2008), pp. 1254--1279


\bibitem{comon2009tensor}
P. Comon, X. Luciani and A. L. De Almeida.
\emph{Tensor decompositions, alternating least squares and other tales. }
Journal of Chemometrics: A Journal of the Chemometrics Society \textbf{23} (2009), pp. 393--405.

\bibitem{de2008tensor}
V. De Silva and L. Lim. 
\emph{Tensor rank and the ill-posedness of the best low-rank approximation problem. }
SIAM Journal on Matrix Analysis and Applications \textbf{30} (2008), pp. 1084--1127.


\bibitem{defant1992tensor}
A. Defant and K. Floret.
\emph{Tensor norms and operator ideals.}
Amsterdam: North-Holland (1993).


\bibitem{diestel2008metric}
J. Diestel, J.H. Fourier and J. Swart.
\emph{The Metric Theory of Tensor Products. Grothendieck's R\'esum\'e Revisited.}
American Mathematical Society (2008).

\bibitem{floret1997natural}
K.~Floret.
\emph{Natural norms on symmetric tensor products of normed spaces.}
Note di Matematica {\bf 17} (1997), pp. 153--188.

\bibitem{friedland2013best}
S. Friedland. 
\emph{Best rank one approximation of real symmetric tensors can be chosen symmetric.}
Frontiers of Mathematics in China \textbf{8} (2013), pp. 19--40.

\bibitem{friedland2014nuclear} 
S. Friedland and L. Lim. 
\emph{Nuclear norm of higher-order tensors.}
Mathematics of Computation {\bf 18} (2018), pp. 1255--1281.

\bibitem{friedland2014number} 
S. Friedland and G. Ottaviani. 
\emph{The number of singular vector tuples and uniqueness of best rank-one approximation of tensors.} 
Foundations of Computational Mathematics \textbf{14} (2014), pp. 1209--1242.

\bibitem{friedland2015low}
S. Friedland and V. Tammali.
\emph{Low-rank approximation of tensors.} 
Numerical Algebra, Matrix Theory, Differential-Algebraic Equations and Control Theory. Springer, Cham (2015), pp. 377--411.

\bibitem{grasedyck2013literature}
L. Grasedyck, D. Kressner and C. Tobler.
\emph{A literature survey of low-rank tensor approximation techniques.}
GAMM--Mitteilungen \textbf{36} (2013), pp. 53--78.

\bibitem{pappas2007norm} A. Pappas, Y. Sarantopoulos and A. Tonge. 
\emph{Norm attaining polynomials}.
Bulletin of the London Mathematical Society {\bf 39} (2007), pp. 255--264.

\bibitem{ryan2002introduction}
R. Ryan. 
\emph{Introduction to Tensor Products of Banach Spaces.}
Springer Monographs in Mathematics (2012).

\bibitem{sidiropoulos2017tensor}
N.~Sidiropoulos, L. De Lathauwer, Xiao Fu, K. Huang, E. Papalexakis and C. Faloutsos. 
\emph{Tensor decomposition for signal processing and machine learning.} 
IEEE Transactions on Signal Processing \textbf{65} (2017), pp. 3551--3582.

\bibitem{zhang2012best}
X. Zhang, C. Ling and L. Qi. 
\emph{The best rank-1 approximation of a symmetric tensor and related spherical optimization problems.}
SIAM Journal on Matrix Analysis and Applications \textbf{33} (2012), pp. 806--821.





\end{thebibliography}
\end{document}